\documentclass[review]{elsarticle}

\pdfoutput=1

\usepackage{lineno,hyperref}
\modulolinenumbers[5]

\journal{Journal}

\graphicspath{{./}{./figs/}}

  \def\clap#1{\hbox to 0pt{\hss#1\hss}}

\usepackage{bm}
\providecommand{\mat}[1]{\bm{#1}}%
\renewcommand{\vec}[1]{\mathbf{#1}}

\providecommand{\mA}{\ensuremath{\mat{A}}}
\providecommand{\mB}{\ensuremath{\mat{B}}}
\providecommand{\mC}{\ensuremath{\mat{C}}}
\providecommand{\mD}{\ensuremath{\mat{D}}}

\providecommand{\mT}{\ensuremath{\mat{T}}}
\providecommand{\mU}{\ensuremath{\mat{U}}}

\providecommand{\mW}{\ensuremath{\mat{W}}}

\providecommand{\va}{\ensuremath{\vec{a}}}
\providecommand{\vb}{\ensuremath{\vec{b}}}

\providecommand{\vq}{\ensuremath{\vec{q}}}
\providecommand{\vr}{\ensuremath{\vec{r}}}

\providecommand{\vu}{\ensuremath{\vec{u}}}
\providecommand{\vv}{\ensuremath{\vec{v}}}
\providecommand{\vw}{\ensuremath{\vec{w}}}
\providecommand{\vx}{\ensuremath{\vec{x}}}
\providecommand{\vy}{\ensuremath{\vec{y}}}

\newcommand{\sS}{\mathcal{S}}

\newcommand{\bmat}[1]{\begin{bmatrix}#1\end{bmatrix}}

\newcommand{\argmin}[1]{\underset{#1}{\mathrm{argmin}}}

\newcommand{\unit}[1]{\text{#1}}
\usepackage{amsmath,amssymb}
\usepackage{subfigure}

\usepackage{kbordermatrix}
\renewcommand{\kbldelim}{[}
\renewcommand{\kbrdelim}{]}

\newcommand{\Vlam}{V_{\text{lam}}}
\newcommand{\Vtur}{V_{\text{tur}}}

\begin{document}

\begin{frontmatter}

\title{Many physical laws are ridge functions}





\author[CSM]{Paul G.~Constantine}
\address[CSM]{Department of Applied Mathematics and Statistics, Colorado School of Mines, 1500 Illinois Street, Golden, CO 80211} \ead{paul.constantine@mines.edu}

\author[SU1]{Zachary del Rosario} 
\address[SU1]{Department of Aeronautics and Astronautics, Stanford University, Durand Building, 496 Lomita Mall, Stanford, CA 94305}

\author[SU2]{Gianluca Iaccarino}
\address[SU2]{Department of Mechanical Engineering and Institute for Computational Mathematical Engineering, Stanford University, Building 500, Stanford, CA 94305}

\begin{abstract}
A ridge function is a function of several variables that is constant along certain directions in its domain. Using classical dimensional analysis, we show that many physical laws are ridge functions; this fact yields insight into the structure of physical laws and motivates further study into ridge functions and their properties. We also connect dimensional analysis to modern subspace-based techniques for dimension reduction, including active subspaces in deterministic approximation and sufficient dimension reduction in statistical regression.
\end{abstract}

\begin{keyword}
active subspaces \sep dimensional analysis \sep dimension reduction \sep sufficient dimension reduction
\end{keyword}

\end{frontmatter}

\linenumbers


In 1969, Harvard Physicist P.~W.~Bridgman~\cite{Bridgman1969} wrote, ``The principal use of dimensional analysis is to deduce from a study of the dimensions of the variables in any physical system certain necessary limitations on the form of any possible relationship between those variables.'' At the time, dimensional analysis was a mature set of tools, and it remains a staple of the science and engineering curriculum because of its ``great generality and mathematical simplicity''~\cite{Bridgman1969}. In this paper, we make Bridgman's ``necessary limitations'' precise by connecting dimensional analysis' fundamental result---the Buckingham Pi Theorem---to a particular low-dimensional structure that arises in modern approximation models.

Today's data deluge motivates researchers across mathematics, statistics, and engineering to pursue exploitable low-dimensional descriptions of complex, high-dimensional systems. Computing advances empower certain structure-exploiting techniques to impact a wide array of important problems. Successes---e.g., compressed sensing in signal processing~\cite{Donoho2006,Candes2006}, neural networks in machine learning~\cite{neuralnets96}, and principal components in data analysis~\cite{pcabook}---abound. In what follows, we review \emph{ridge functions}~\cite{pinkus2015}, which exhibit a particular type of low-dimensional structure, and we show how that structure manifests in physical laws.

Let $\vx\in\mathbb{R}^m$ be a vector of continuous parameters; a ridge function $f:\mathbb{R}^m\rightarrow\mathbb{R}$ takes the form
\begin{equation}
\label{eq:ridgedef}
f(\vx) \;=\; g(\mA^T\vx),
\end{equation}
where $\mA\in\mathbb{R}^{m\times n}$ is a constant matrix with $n<m$, and $g:\mathbb{R}^{n}\rightarrow\mathbb{R}$ is a scalar-valued function of $n$ variables. Although $f$ is nominally a function of $m$ variables, it is constant along all directions orthogonal to $\mA$'s columns. To see this, let $\vx\in\mathbb{R}^m$ and $\vy=\vx+\vu\in\mathbb{R}^m$ with $\vu$ orthogonal to $\mA$'s columns, i.e., $\mA^T\vu=0$. Then
\begin{equation}
\label{eq:constant}
f(\vy) 
\;=\; g(\mA^T(\vx + \vu)) 
\;=\; g(\mA^T\vx)
\;=\; f(\vx).
\end{equation}
Ridge functions appear in multivariate Fourier transforms, plane waves in partial differential equations, and statistical models such as projection pursuit regression and neural networks; see~\cite[Chapter 1]{pinkus2015} for a comprehensive introduction. Ridge functions have recently become an object of study in approximation theory~\cite{Donoho2001}, and computational scientists have proposed methods for estimating their properties (e.g., the columns of $\mA$ and the form of $g$) from point evaluations $f(\vx)$~\cite{Fornasier2012,Cohen11}. However, scientists and engineers outside of mathematical sciences have paid less attention to ridge functions than other useful forms of low-dimensional structure. Many natural signals are sparse, and many real world data sets contain colinear factors. But whether ridge structures are pervasive in natural phenomena remains an open question. 

We answer this question affirmatively by showing that \emph{many physical laws are ridge functions}. This conclusion is a corollary of classical dimensional analysis. To show this result, we first review classical dimensional analysis from a linear algebra perspective. 

\section{Dimensional analysis}
\label{sec:da}

Several physics and engineering textbooks describe classical dimensional analysis. Barenblatt~\cite{Barenblatt1996} provides a thorough overview in the context of scaling and self-similarity, while Ronin~\cite{sonin2001} is more concise. However, the presentation by Calvetti and Somersalo~\cite{calvetti15}, which ties dimensional analysis to linear algebra, is most appropriate for our purpose; what follows is similar to their treatment.  

We assume a chosen measurement system has $k$ fundamental units. For example, a mechanical system may have the $k=3$ fundamental units of time (seconds, $\unit{s}$), length (meters, $\unit{m}$), and mass (kilograms, $\unit{kg}$). More generally, for a system in SI units, $k$ is at most 7---the base units. All measured quantities in the system have units that are products of powers of the fundamental units; for example, velocity has units of length per time, $\unit{m}\cdot\unit{s}^{-1}$. 

Define the \emph{dimension function} of a quantity $q$, denoted $[q]$, to be a function that returns the units of $q$; if $q$ is dimensionless, then $[q]=1$. Define the \emph{dimension vector} of a quantity $q$, denoted $\vv(q)$, to be a function that returns the $k$ exponents of $[q]$ with respect to the $k$ fundamental units; if $q$ is dimensionless, then $\vv(q)$ is a $k$-vector of zeros. For example, in a system with fundamental units $\unit{m}$, $\unit{s}$, and $\unit{kg}$, if $q$ is velocity, then $[q]=\unit{m}^1\cdot\unit{s}^{-1}\cdot\unit{kg}^0$ and $\vv(q)=[1,-1,0]^T$; the order of the units does not matter. 

Barenblatt~\cite[Section 1.1.5]{Barenblatt1996} states that quantities $q_1,\dots,q_m$ ``have \emph{independent dimensions} if none of these quantities has a dimension function that can be represented in terms of a product of powers of the dimensions of the remaining quantities.'' This is equivalent to linear independence of the associated dimension vectors $\{\vv(q_1),\dots,\vv(q_m)\}$. If quantities $q_1,\dots,q_m$ have independent dimensions with respect to $k$ fundamental units, then $m\leq k$. We say that the quantities' independent dimensions are \emph{complete} if $m=k$. We can express the exponents for a derived dimension as a linear system of equations. Let $q_1,\dots,q_m$ contain quantities with complete and independent dimensions, and define the $k\times m$ matrix
\begin{equation}
\label{eq:D}
\mD \;=\; \bmat{\vv(q_1) &\cdots & \vv(q_m)}.
\end{equation}
By independence, $\mD$ has rank $k$. Let $p$ be a dimensional quantity with derived units $[p]$. Then $[p]$ can be written as products of powers of $[q_1],\dots,[q_m]$,
\begin{equation}
\label{eq:dimeq}
[p] \;=\; [q_1]^{w_1}\cdots[q_m]^{w_m}.
\end{equation}
The powers $w_1,\dots,w_m$ satisfy the linear system of equations
\begin{equation}
\label{eq:linsys}
\mD\,\vw = \vv(p), \qquad
\vw = \bmat{w_1\\ \vdots \\ w_m}. 
\end{equation}
Given the solution $\vw$ of \eqref{eq:linsys}, we can define a quantity $p'$ with the same units as $p$ (i.e., $[p]=[p']$) as
\begin{equation}
\label{eq:logtrans}
\begin{aligned}
p' &= q_1^{w_1}\cdots q_m^{w_m}\\
&= \exp\left(\log\left(
q_1^{w_1}\cdots q_m^{w_m}
\right)\right)\\
&= \exp\left(\sum_{i=1}^m
w_i\,\log(q_i)\right)\\
&= \exp\left(\vw^T\log(\vq)\right),
\end{aligned}
\end{equation}
where $\vq=[q_1,\dots,q_m]^T$, and the log of a vector returns the log of each component. There is some controversy over whether the logarithm of a physical quantity makes physical sense~\cite{molyneux91}. We sidestep this discussion by noting that $\exp(\vw^T\log(\vq))$ is merely a formal expression of products of powers of physical quantities. There is no need to interpret the units of the logarithm of a physical quantity. 

\subsection{Nondimensionalization}

Assume we have a system with $m+1$ dimensional quantities, $q$ and $\vq=[q_1,\dots,q_m]^T$, whose dimensions are derived from a set of $k$ fundamental units, and $m>k$. Without loss of generality, assume that $q$ is the quantity of interest with units $[q]$, and we assume $[q]\not=1$ (i.e., that $q$ is not dimensionless). We assume that $\mD$, defined as in \eqref{eq:D}, has rank $k$, which is equivalent to assuming that there is a complete set of dimensions among $[q_1],\dots,[q_m]$. We construct a dimensionless quantity of interest $\pi$ as 
\begin{equation}
\label{eq:qoi}
\pi \;=\; \pi(q,\,\vq) \;=\; q\,\exp(-\vw^T\log(\vq)),
\end{equation}
where the exponents $\vw$ satisfy the linear system $\mD\vw=\vv(q)$. The solution $\vw$ is not unique, since $\mD$ has a nontrivial null space. 

Let $\mW=[\vw_1,\dots,\vw_n]\in\mathbb{R}^{m\times n}$ be a matrix whose columns contain a basis for the null space of $\mD$. In other words,
\begin{equation}
\label{eq:null0}
\mD\mW \;=\; \mathbf{0}_{k\times n},
\end{equation}
where $\mathbf{0}_{k\times n}$ is an $k$-by-$n$ matrix of zeros. The completeness assumption implies that $n=m-k$. Note that the basis for the null space is not unique, which is a challenge in classical analysis.
Calvetti and Somersalo~\cite[Chapter 4]{calvetti2013} offer a recipe for computing $\mW$ with rational elements via Gaussian elimination; this is consistent with the physically intuitive construction of many classical nondimensional quantities such as the Reynolds number. However, one must still choose the pivot columns in the Gaussian elimination.

The Buckingham Pi Theorem~\cite[Chapter 1.2]{Barenblatt1996} states that any physical law can be expressed as a relationship between the dimensionless quantity of interest $\pi$ and the $n=m-k$ dimensionless quantities. Similar to \eqref{eq:logtrans}, we can formally express the dimensionless parameters $\pi_i$ as
\begin{equation}
\label{eq:dimless}
\pi_i \;=\; \pi_i(\vq) \;=\; \exp(\vw_i^T\log(\vq)),\qquad i=1,\dots,n.
\end{equation}
The dimensionless parameters depend on the choice of basis vectors $\vw_i$. We seek a function $f:\mathbb{R}^{n}\rightarrow\mathbb{R}$ that models the relationship between $\pi$ and $\pi_1,\dots,\pi_n$,
\begin{equation}
\label{eq:phi}
\pi \;=\; f(\pi_1,\dots,\pi_n).
\end{equation}
Expressing the physical law in dimensionless quantities has several advantages. First, there are typically fewer dimensionless quantities than measured dimensional quantities, which allows one to construct $f$ with many fewer experiments than one would need to build a relationship among dimensional quantities; several classical examples showcase this advantage~\cite[Chapter 1]{Barenblatt1996}. Second, dimensionless quantities do not change if units are scaled, which allows one to devise small scale experiments that reveal a scale-invariant relationship. Third, when all quantities are dimensionless, any mathematical relationship will satisfy \emph{dimensional homogeneity}, which is a physical requirement that models only sum quantities with the same dimension.

\section{Physical laws are ridge functions}
\label{sec:pi}

We exploit the nondimensionalized statement of the physical law to show that the corresponding dimensional physical law can be expressed as a ridge function of dimensional quantities by combining \eqref{eq:qoi}, \eqref{eq:phi}, and \eqref{eq:dimless} as follows:
\begin{equation}
\label{eq:rd0}
\begin{aligned}
q\,\exp(-\vw^T\log(\vq))
&= \pi \\
&= f(\pi_1,\dots,\pi_n)\\
&= f\left(
\exp(\vw_1^T\log(\vq)),\dots,\exp(\vw_n^T\log(\vq))
\right).
\end{aligned}
\end{equation}
We rewrite this expression as
\begin{equation}
\label{eq:rd3}
\begin{aligned}
q &= \exp(\vw^T\log(\vq))\,f\left(
\exp(\vw_1^T\log(\vq)),\dots,\exp(\vw_n^T\log(\vq))
\right)\\
&= h(\vw^T\log(\vq),\vw_1^T\log(\vq),\dots,\vw_n^T\log(\vq))\\
&= h(\mA^T\vx),
\end{aligned}
\end{equation}
where $h:\mathbb{R}^{n+1}\rightarrow\mathbb{R}$, and the variables $\vx=\log(\vq)$ are the logs of the dimensional quantities. The matrix $\mA$ contains the vectors computed in \eqref{eq:qoi} and \eqref{eq:dimless},
\begin{equation}
\label{eq:A}
\mA \;=\; 
\bmat{\vw & \mW} \;=\;
\bmat{\vw&\vw_1&\cdots &\vw_n}\in\mathbb{R}^{m\times (n+1)}.
\end{equation}
The form \eqref{eq:rd3} is a ridge function in $\vx$, which justifies our thesis; compare to \eqref{eq:ridgedef}. We call the column space of $\mA$ the \emph{dimensional analysis subspace}. 

Several remarks are in order. First, \eqref{eq:rd3} reveals $h$'s dependence on its first coordinate. If one tries to fit $h(y_0,y_1,\dots,y_n)$ from measured data---as in semi-empirical modeling---then she should pursue a function of the form $h(y_0,y_1,\dots,y_n)=\exp(y_0)\,g(y_1,\dots,y_n)$, where $g:\mathbb{R}^n\rightarrow\mathbb{R}$. In other words, the first input of $h$ merely scales another function of the remaining variables. 

Second, writing the physical law as $q=h(\mA^T\vx)$ as in \eqref{eq:rd3} uses a ridge function of the logs of the physical quantities. In dimensional analysis, some contend that the log of a physical quantity is not physically meaningful. However, there is no issue taking the log of numbers; measured data is often plotted on a log scale, which is equivalent. To construct $h$ from measured data, one must compute the logs of the measured numbers; such fitting is a computational exercise that ignores the quantities' units. 

To evaluate a fitted semi-empirical model, we compute $q$ given $\vq$ as
\begin{equation}
\begin{aligned}
q &= h(\vw^T\log(\vq),\vw_1^T\log(\vq),\dots,\vw_n^T\log(\vq))\\
&= \exp(\vw^T\log(\vq)) \,g\big(\log(\exp(\vw_1^T\log(\vq))),\dots, \log(\exp(\vw_n^T\log(\vq)))\big)\\
&= \exp(\vw^T\log(\vq))\,g\big(\log(\pi_1),\dots,\log(\pi_n)\big)
\end{aligned}
\end{equation}
In $g\left(\log(\pi_1),\dots,\log(\pi_n)\right)$, the logs take nondimensional quantities, and $g$ returns a nondimensional quantity. By construction, the term $\exp(\vw^T\log(\vq))$ has the same dimension as $q$, so dimensional homogeneity is satisfied. Thus, the ridge function form of the physical law \eqref{eq:rd3} does not violate dimensional homogeneity. 

Third, the columns of $\mA$ are linearly independent by construction. The first column is not in the null space of $\mD$ (see \eqref{eq:qoi}), and the remaining columns form a basis for the null space of $\mD$ (see \eqref{eq:null0}). So $\mA$ has full column rank. Then the quantity of interest $q$ is invariant to changes in the inputs $\vx$ that live in the null space of $\mA^T$; see \eqref{eq:constant}. 

\section{Relationships to other subspaces}
\label{sec:dr}

We have shown that physical laws are ridge functions as a consequence of classical dimensional analysis. This observation connects physical modeling to two modern analysis techniques: \emph{active subspaces} in deterministic approximation and \emph{sufficient dimension reduction} in statistical regression. We note that recent statistics literature has explored the importance of dimensional analysis for statistical analyses, e.g., design of experiments~\cite{Albrecth2013} and regression analysis~\cite{Shen2014}. These works implicitly exploit the ridge-like structure in the physical laws; in what follows, we make these connections explicit.

\subsection{Active subspaces}

The active subspace of a given function is defined by a set of important directions in the function's domain. More precisely, let $f(\vx)$ be a differentiable function from $\mathbb{R}^m$ to $\mathbb{R}$, and let $p(\vx)$ be a bounded probability density function on $\mathbb{R}^m$. Define the $m\times m$ symmetric and positive semidefinite matrix $\mC$ as
\begin{equation}
\label{eq:Cmat}
\mC \;=\; \int \nabla f(\vx)\,\nabla f(\vx)^T\,p(\vx)\,d\vx,
\end{equation}
where $\nabla f(\vx)$ is the gradient of $f$. The matrix $\mC$ admits a real eigenvalue decomposition $\mC=\mU\Lambda\mU^T$, where $\mU$ is the orthogonal matrix of eigenvectors, and $\Lambda$ is the diagonal matrix of non-negative eigenvalues denoted $\lambda_1,\dots,\lambda_m$ ordered from largest to smallest. Assume that $\lambda_k>\lambda_{k+1}$ for some $k<m$. Then the active subspace of $f$ is the span of the first $k$ eigenvectors. Our recent work has developed computational procedures for first estimating the active subspace and then exploiting it to enable calculations that are otherwise prohibitively expensive when the number $m$ of components in $\vx$ is large---e.g., approximation, optimization, and integration~\cite{asbook,constantine2014active}. 

If $f$ is a ridge function as in \eqref{eq:ridgedef}, then $f$'s active subspace is related to the $m\times n$ matrix $\mA$. First, observe that $\nabla f(\vx) = \mA\,\nabla g(\mA^T\vx)$, where $\nabla g$ is the gradient of $g$ with respect to its arguments. Then $\mC = \mA\mT\mA^T$, where
\begin{equation}
\mT\;=\;
\int \nabla g(\mA^T\vx)\,\nabla g(\mA^T\vx)^T\,p(\vx)\,d\vx.
\end{equation}
The symmetric positive semidefinite matrix $\mT$ has size $n\times n$. The form of $\mC$ implies two facts: (i) $\mC$ has rank at most $n$, and (ii) the invariant subspaces of $\mC$ up to dimension $n$ are subspaces of $\mA$'s column space. Moreover, if $\mA$'s columns are orthogonal, then we can compute the first $n$ of $\mC$'s eigenpairs by computing $\mT$'s eigenpairs, which offers a computational advantage.

Since a physical law is a ridge function (see \eqref{eq:rd3}), the second fact implies that the active subspace of $q$ is a subspace of the dimensional analysis subspace. Additionally, if the functional form is transformed so that $\mA$ has orthogonal columns (e.g., via a QR factorization), then one may estimate $q$'s active subspace with less effort by exploiting the connection to dimensional analysis. 

\subsection{Sufficient dimension reduction}

The tools associated with active subspaces apply to deterministic approximation problems. In a physical experiment, where measurements are assumed to contain random noise, statistical regression may be a more appropriate tool. Suppose that $N$ independent experiments each produce measurements $(\vx_i,y_i)$; in the regression context, $\vx_i$ is the $i$th sample of the predictors and $y_{i}$ is the associated response with $i=1,\dots,N$. These quantities are related by the regression model,
\begin{equation}
\label{eq:regress}
y_{i} \;=\; f(\vx_i) + \varepsilon_i,
\end{equation}
where $\varepsilon_i$ are independent random variables. Let $F_{y|\vx}(\cdot)$ be the cumulative distribution function of the random variable $y$ conditioned on $\vx$. Suppose $\mB\in\mathbb{R}^{m\times k}$ is such that 
\begin{equation}
F_{y|\vx}(a) \;=\; F_{y|\mB^T\vx}(a), \qquad a\in\mathbb{R}.
\end{equation} 
In words, the information about $y$ given $\vx$ is the same as the information about $y$ given linear combinations of the predictors $\mB^T\vx$. When this happens, the column space of $\mB$ is called a \emph{dimension reduction subspace}~\cite{cook2009regression}. There is a large body of statistics literature that describes methods for estimating the dimension reduction subspace given samples of predictor/response pairs; see Cook~\cite{cook2009regression} for a comprehensive review. These methods fall under the category of \emph{sufficient dimension reduction}, since the dimension reduction subspace is sufficient to statistically characterize the regression.

The ridge function structure in the physical law \eqref{eq:rd3} implies that the dimensional analysis subspace is a dimension reduction subspace, where logs of the physical quantities $\vx$ are the predictors, and the quantity of interest $q$ is the response. Note that the dimensional analysis subspace may not be a \emph{minimal} dimension reduction subspace---i.e., a dimension reduction subspace with the smallest dimension. This connection may lead to new or improved sufficient dimension reduction methods that incorporate dimensional analysis. 

\section{Example: viscous pipe flow}

To demonstrate the relationship between the dimensional analysis subspace and the active subspace, we consider the classical example of viscous flow through a pipe. The system's three fundamental units ($k=3$) are kilograms (kg), meters (m), and seconds (s). The physical quantities include the fluid's bulk velocity $V$, density $\rho$, and viscosity $\mu$; the pipe's diameter $D$ and characteristic wall roughness $\varepsilon$; and the pressure gradient $\frac{\Delta P}{L}$. We treat $V$ as the quantity of interest. 


\subsection{Dimensional analysis}
The matrix $\mD$ from \eqref{eq:D} encodes the units; for this system $\mD$ is
\begin{equation}
\kbordermatrix{ & \rho & \mu & D & \varepsilon & \Delta P/L \cr
\text{kg} & 1 & 1 & 0 & 0 & 1\cr
\text{m} & -3 & -1 & 1 & 1 & -2\cr 
\text{s} & 0 & -1 & 0 & 0 & -2 }
\end{equation}
The vector $\vw$ from \eqref{eq:qoi} that nondimensionalizes the velocity is $[-2,1,0,0,1]^T$. The matrix $\mW$ whose columns span the null space of $\mD$ is 
\begin{equation}
\mW \;=\; \bmat{0 & 1\\ 0 & -2\\ -1 & 3\\ 1 & 0\\ 0 & 1}.
\end{equation} 
Thus, the dimensional analysis subspace is the span of $\vw$ and $\mW$'s two columns.

\subsection{Active subspace}
To estimate the active subspace for this system, we wrote a MATLAB code to compute velocity as a function the other physical quantities. The code relies on well established theories for this system (\cite[Chapter 6]{White2011} and \cite{Moody1944,Colebrook1937}); see supporting information for more details. We compute gradients with a first order finite difference approximation. We compare the active subspace across two parameter regimes: one corresponding to laminar flow and the other corresponding to turbulent flow. We characterize the regimes by a range on each of the physical quantities; details are in the supporting information. In each case, the probability density function $p(\vx)$ from \eqref{eq:Cmat} is a uniform density on the space defined by the parameter ranges. We estimate the integrals defining $\mC$ from \eqref{eq:Cmat} using a tensor product Gauss-Legendre quadrature rule with 11 points in each dimension (161051 points in five dimensions), which was sufficient for 10 digits of accuracy in the eigenvalues. 

\begin{figure}[ht]
\centering
\begin{subfigure}{}
\includegraphics[width=.43\textwidth]{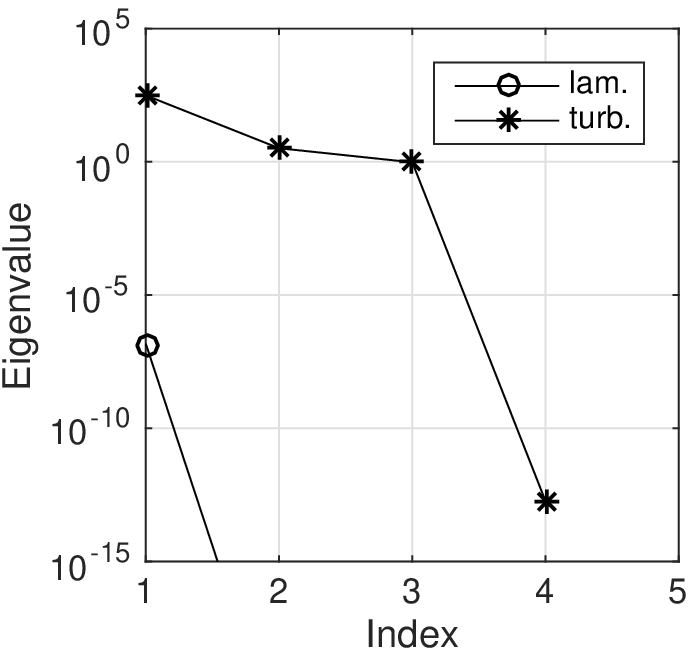}
\label{fig:eigs}
\end{subfigure}
\begin{subfigure}{}
\includegraphics[width=.43\textwidth]{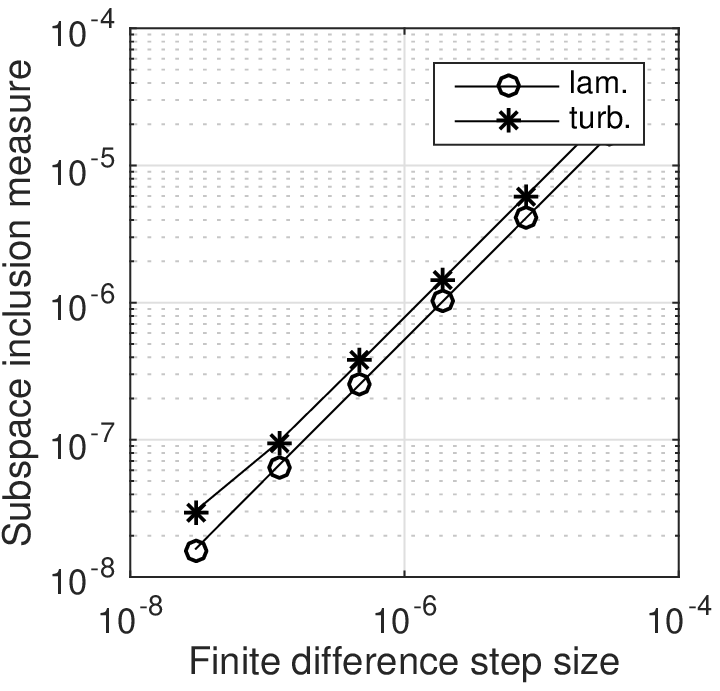}
\label{fig:errs}
\end{subfigure}
\caption{The left figure shows the eigenvalues of $\mC$ from \eqref{eq:Cmat} for the laminar and turbulent regimes. The laminar regime has one nonzero eigenvalue (to numerical precision), and the turbulent regime has three (to numerical precision). The right figure measures the inclusion of the one-dimensional (laminar) and three-dimensional (turbulent) active subspaces in the three-dimensional dimensional analysis subspace as a function of the finite difference step size.}
\label{fig:e}
\end{figure}

For laminar flow, the velocity is equal to a product of powers of the other quantities. Therefore, we expect only one non-zero eigenvalue for $\mC$ from \eqref{eq:Cmat}, which indicates a one-dimensional active subspace. For turbulent flow, the relationship is more complex, so we expect up to three non-zero eigenvalues for $\mC$. The left subfigure in Figure \ref{fig:e} shows the eigenvalues of $\mC$ for these two cases. Indeed, there is one relatively large eigenvalue in the laminar case and three relatively large eigenvalues in the turbulent case. The fourth eigenvalue in the turbulent regime is roughly $10^{-13}$, which is within the numerical accuracy of the integrals and gradients. The right subfigure in Figure \ref{fig:e} shows the amount by which the active subspace---one-dimensional in the laminar case and three-dimensional in the turbulent case---is not a subset of the three-dimensional dimensional analysis subspace as a function of the finite difference step size; details on this measurement are in the supporting information. Note the first order convergence of this metric toward zero as the finite difference step size decreases. This provides strong numerical evidence that the active subspace is subset of the dimensional analysis subspace in both cases as the theory predicts. 

\section{Conclusions}
\label{sec:conc}

We have shown that classical dimensional analysis implies that many physical laws are ridge functions. This fact motivates further study into ridge functions---both analytical and computational. The result is a statement about the general structure of physical laws. We expect there are many ways a modeler may exploit this structure, e.g., for building semi-empirical models from data or finding insights into invariance properties of the physical system. We also connect the ridge function structure to modern subspace-based dimension reduction ideas in approximation and statistical regression. We hope that this explicit connection motivates modelers to explore these techniques for finding low-dimensional parameterizations of complex, highly parameterized models. 

The dimension reduction enabled by the dimensional analysis subspace is naturally limited by the number of fundamental units. The SI units contain seven base units. Therefore, the dimension of the subspace in which the physical law is invariant (i.e., the dimension of the complement of the dimensional analysis subspace) is at most six in any physical system with SI units. For many systems, reducing the number of input parameters by six will be remarkably beneficial---potentially enabling studies and experiments not otherwise feasible. However, there may be other systems where reducing the dimension by six still yields an intractable reduced model. In such cases, the modeler may explore the subspace-based dimension reduction techniques for potentially greater reduction. The connections we have established between these techniques and the dimensional analysis subspace aid in efficient implementation and interpretation of results.

\section{Supporting information}

\subsection{Viscous pipe flow}

To study the relationship between the active subspace and the dimensional analysis subspace, we consider the classical example of a straight pipe with circular cross-section and rough walls filled with a viscous fluid. A pressure gradient is applied which drives axial flow. The system's three fundamental units $(k=3)$ are kilograms (kg), meters (m), and seconds (s). The physical quantities include the fluid's bulk velocity $V$, density $\rho$, and viscosity $\mu$; the pipe's diameter $D$ and characteristic wall roughness $\varepsilon$; and the pressure gradient $\frac{\Delta P}{L}$. We treat velocity $V$ as the quantity of interest, noting that it is equal to the volumetric flow rate through the pipe divided by the cross-sectional area.

\begin{figure}[h]
\centering
\includegraphics[width=.71\textwidth]{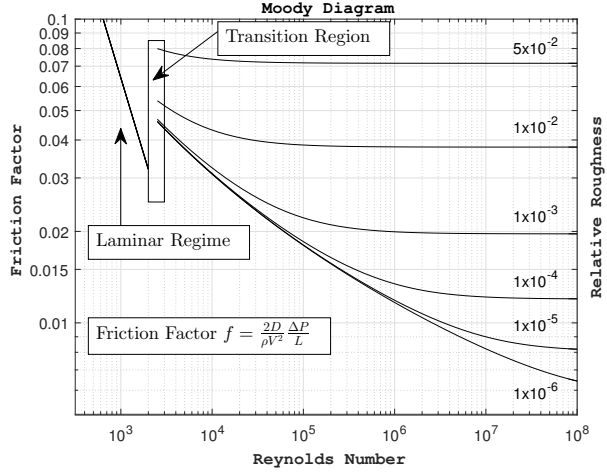}
\caption{The Moody Diagram plots the friction factor (dimensionless pressure loss) against the Reynolds number and relative roughness. Transition from laminar flow (governed by the Poiseuille relation) to turbulent flow (modeled by the Colebrook equation) is assumed to occur at a critical Reynolds number $Re_c \approx 3\times10^3$.}
\label{fig:moody}
\end{figure}

These quantities are implicitly related by the Moody Diagram (Figure \ref{fig:moody}), which plots the friction factor $f$ defined by
\begin{equation}
f \;=\; \frac{\Delta P}{L}\frac{D}{\frac{1}{2}\rho V^2},
\end{equation}
against the Reynolds number $\frac{\rho V D}{\mu}$ and relative roughness $\frac{\varepsilon}{D}$~\cite{Moody1944}. Below a critical Reynolds number, taken to be $Re_c=3\times10^3$, the friction factor satisfies the Poiseuille relation~\cite[Chapter 6]{White2011},
\begin{equation}
\label{eq:poiseuille}
f \;=\; \frac{64}{Re}, 
\end{equation}
For $Re>Re_c$, the Colebrook equation~\cite{Colebrook1937} implicitly defines the relationship between friction factor and the other quantities,
\begin{equation}
\label{eq:colebrook}
\frac{1}{\sqrt{f}} \;=\; -2.0\log_{10}\left(\frac{1}{3.7}\frac{\varepsilon}{D} + \frac{2.51}{Re\sqrt{f}}\right), 
\end{equation}
This relationship is valid through transition to full turbulence.

\subsection{Bulk velocity as the quantity of interest}

Substituting dimensional quantities and solving for $V$ in \eqref{eq:poiseuille} yields an expression for bulk velocity in laminar flow, denoted $\Vlam$, 
\begin{equation}
\label{eq:laminar_velocity}
\Vlam \;=\; \frac{\Delta P}{L}\frac{D^2}{32\mu}.
\end{equation}
Note that the expression on the right hand side is exactly a product of powers of the remaining dimensional quantities. Thus, we can write the right hand side in a form similar to [\textbf{6}] in the main manuscript, which shows that $\Vlam$ is a ridge function of one linear combination of the logs of the dimensional quantities. And we therefore expect laminar velocity to have a one-dimensional active subspace, despite the fact that the more generic dimensional analysis subspace is three-dimensional. We verify this in the numerical experiment represented by Figure 1 in the main manuscript. 

Substituting dimensional quantities in \eqref{eq:colebrook}, and noting that $V$ cancels within logarithm, reveals the explicit relation for bulk velocity in turbulent flow $\Vtur$,
\begin{equation}
\label{eq:turbulent_velocity}
\Vtur \;=\; -2.0\sqrt{\frac{\Delta P}{L}\frac{2D}{\rho}} \log_{10}\left(%
    \frac{1}{3.7}\frac{\varepsilon}{D} + 2.51\frac{\mu}{D^{3/2}}%
    \sqrt{\frac{L}{\Delta P}\frac{1}{2\rho}}\right).
\end{equation}
Note that the right hand side of \eqref{eq:turbulent_velocity} is more complicated than the right hand side of \eqref{eq:laminar_velocity} due to the logarithm term; i.e., it is more than a product of powers. Therefore, we expect that the active subspace for turbulent bulk velocity has dimension greater than one but not more than three---since it is a subspace of the three-dimensional dimensional analysis subspace. Figure 1 from the main manuscript numerically verifies this observation.

Given inputs $\rho$, $D$, $\mu$, $\frac{\Delta P}{L}$, and $\varepsilon$, we obtain $V$ by choosing between $\Vlam$ and $\Vtur$. To make this choice, we compute $Re$ based on $\Vtur$ and set $V$ to $\Vtur$ if this value exceeds $Re_{c}$. Otherwise, we set $V$ to be $\Vlam$. We wrote a MATLAB script to reproduce these relationships, and we treat the script as a virtual laboratory that we use to verify that the active subspaces---one for turbulent flow and one for laminar flow---satisfy the theoretical relationship to the dimensional analysis subspace. 

\subsection{Parameter ranges}

To define the active subspace for bulk velocity, we need a density function on the logs of the input quantities. We choose the density function to be a uniform, constant density on a five-dimensional hyperrectangle, defined by ranges on the input quantities, and zero elsewhere. We choose the input ranges to produce flow that is essentially laminar or essentially turbulent---depending on the associated Reynolds number.  

\begin{table}[h]
\centering
\caption{Parameter bounds for the laminar flow case.}
\label{tab:bounds_laminar}
\begin{tabular}{lllll}
fluid density     & $\rho$ & $1.0\times10^{-1}$ & $1.4\times10^{-1}$ & $\unit{kg}/\unit{m}^3$ \\
fluid viscosity   & $\mu$  & $1.0\times10^{-6}$   & $1.0\times10^{-5}$   & $\unit{kg} / (\unit{m}\unit{s})$ \\
pipe diameter     & $D$    & $1.0\times10^{-1}$   & $1.0\times10^{+0}$   & $\unit{m}$ \\
pipe roughness    & $\varepsilon$ & $1.0\times10^{-3}$ & $1.0\times10^{-1}$ & $\unit{m}$ \\
pressure gradient & $\frac{\Delta P}{L}$ & $1.0\times10^{-9}$ & $1.0\times10^{-7}$ & $\unit{kg} / (\unit{m}\unit{s})^2$
\end{tabular}
\end{table}

\begin{table}[h]
\centering
\caption{Parameter bounds for the turbulent flow case.}
\label{tab:bounds_turbulent}
\begin{tabular}{lllll}
fluid density     & $\rho$ & $1.0\times10^{-1}$ & $1.4\times10^{-1}$ & $\unit{kg}/\unit{m}^3$ \\
fluid viscosity   & $\mu$  & $1.0\times10^{-6}$   & $1.0\times10^{-5}$   & $\unit{kg} / (\unit{m}\unit{s})$ \\
pipe diameter     & $D$    & $1.0\times10^{-1}$   & $1.0\times10^{+0}$   & $\unit{m}$ \\
pipe roughness    & $\varepsilon$ & $1.0\times10^{-3}$ & $1.0\times10^{-1}$ & $\unit{m}$ \\
pressure gradient & $\frac{\Delta P}{L}$ & $1.0\times10^{-1}$ & $1.0\times10^{+1}$ & $\unit{kg} / (\unit{m}\unit{s})^2$
\end{tabular}
\end{table}

Table \ref{tab:bounds_laminar} shows the parameter ranges that result in laminar flow, and Table \ref{tab:bounds_turbulent} shows the parameter ranges that result in essentially turbulent flow. Approximately 98\% of the Gaussian quadrature points used to estimate integrals for the turbulent case produce turbulent flow cases, i.e., a Reynolds number above the critical threshold. 

\subsection{Subspace inclusion}

For the velocity model, we compute a basis for the dimensional analysis subspace, and we use numerical integration and numerical differentiation to estimate the matrix
\begin{equation}
\mC \;=\; \int \nabla f(\vx)\,\nabla f(\vx)^T\,p(\vx)\,d\vx,
\end{equation}
where $f$ represents the bulk velocity, $\vx$ represents the logs of the remaining dimensional quantities, and $p(\vx)$ is the density function for one of the two flow cases, laminar or turbulent. We approximate the active subspace using the numerical estimates of the eigenpairs of the numerical estimate of $\mC$. To show that, in both cases, the active subspace is a subspace of the dimensional analysis subspace, we use the following numerical test. Consider two subspaces, $\sS_1\subset\mathbb{R}^n$ and $\sS_2\subset\mathbb{R}^m$, with respective bases $\mB_1=[\vb_{1,1},\dots,\vb_{1,n}]$ and $\mB_2=[\vb_{2,1},\dots,\vb_{2,m}]$, where $n<m$. To check if $\sS_1$ is a subspace of $\sS_2$, it is sufficient to check if each column of $\mB_1$ can be represented as a linear combination of the columns of $\mB_2$. Define the residuals
\begin{equation}
\vr_i \;=\; \mB_2 \va_i^* - \vb_{1,i}, \quad i=1,\dots,n,
\end{equation}
where $\va_i^*$ is the minimizer
\begin{equation}
\va_i^* \;=\; \argmin{\va\in\mathbb{R}^m}\, \frac{1}{2}\,\|\mB_2 \va - \vb_{1,i} \|_2^2.
\end{equation}
Define the total residual norm $r^2$ as
\begin{equation}
r^2 \;=\; \sum_{i=1}^n \|\vr_i\|^2_2.
\end{equation}
If $r^2=0$, then $\sS_1$ is a subspace of $\sS_2$. For our numerical example, the errors due to Gaussian quadrature are negligible; we have used enough points to ensure 10 digits of accuracy in all quantities. However, the errors due to numerical differentiation is not negligible. Our numerical test shows that $r^2$ converges to zero like $O(h)$, where $h$ is the finite difference step size, as expected for a first order finite difference approximation. This provides evidence that the active subspace is a subspace of the dimensional analysis subspace, for both flow cases, as numerical errors decrease.

\section*{Acknowledgments}

This material is based on work supported by Department of Defense, Defense Advanced Research Project Agency's program Enabling Quantification of Uncertainty in Physical Systems. The second author's work is supported by the National Science Foundation Graduate Research Fellowship under Grant No. DGE-114747.

\section*{References}

\bibliography{pi-subspaces}

\end{document}